\documentclass[a4paper,10pt]{article}
\usepackage{a4wide}
\usepackage{amsmath}
\usepackage{amssymb}
\usepackage{amsthm}
\usepackage{latexsym}
\usepackage{graphicx}
\usepackage[english]{babel}
\usepackage{makeidx}

\newtheorem{theorem}{Theorem}[section]
\newtheorem{lemma}[theorem]{Lemma}
\newtheorem{corollary}[theorem]{Corollary}
\newtheorem{proposition}[theorem]{Proposition}

\theoremstyle{definition}
\newtheorem{definition}[theorem]{Definition}
\newtheorem{example}[theorem]{Example}

\begin{document}
\title{The type of the base ring associated to a transversal polymatroid
\thanks{This paper was supported by CIDAGT 249-Cex-M3-249-2006}}

\author{Alin \c{S}tefan\\
 "Petroleum and Gas" University of Ploie\c sti, Romania}
\date{}
\maketitle
\begin{abstract}
In this paper we determine the $type,$  $a-invariant$, 
$Hilbert \ function$ 
and $Hilbert \ series$ of the base ring of some transversal 
polymatroid such that its  associated cone has dimension $n$ and $n+1$ facets.  
\end{abstract}

\begin{quotation}


\end{quotation}

\section{Introduction}
In this paper we determine the facets of the polyhedral cone generated
 by the exponent set of the monomials defining the base ring associated
 to some transversal polymatroid. We need the description of these facets
 to find the canonical module of the base ring which is
 expressed in terms of the relative interior of the cone. This would allow
 us to compute the $a$-invariant of those base rings. The results
 presented were discovered by extensive computer algebra experiments
performed with {\it{Normaliz}} \ \cite{BK}.

\section{Preliminaries}
\[\]
Let $n\in \mathbb{N},$ \ $n \geq 3,$ $\sigma \in S_{n},$ \
$\sigma=(1,2,\ldots, n)$ \ the cycle of length $n,$ \ $[n]:=\{1,
2,\ldots, n\}$ \ and $\{e_{i}\}_{1 \leq i \leq n}$ \ be the
canonical base of $\mathbb{R}^{n}.$  For a vector $x\in
 \mathbb{R}^{n}$, $x=(x_{1},\ldots,x_{n})$,  we will denote  $| \ x \ | := x_{1}+
 \ldots + x_{n}.$ If $x^{a}$ is a monomial in $K[x_{1}, \ldots, x_{n}]$ we
 set $\log (x^{a})=a$. Given a set $A$ of monomials, the $log \ set \ of \
 A,$ denoted $\log (A),$ consists of all  $\log (x^{a})$ with $x^{a}\in
 A.$  
We consider the following set of
integer vectors of $\mathbb{N}^{n}$:

\[\begin{array}{cccccccccccccccccc}
     \downarrow i^{th} column  \ \ \ \ \ \ \ \ \ \ \ \ \ \ \ 
     \ \ \ \ \ \ \ \ \ 
  \end{array}
\]

$\nu^{j}_{\sigma^{0}[i] } := \left( \begin{array}{cccccccccccc}
    -j, & -j, & \ldots & ,-j, & (n-j), & \ldots & ,(n-j)  \\
   \end{array}
   \right),
$

\[\begin{array}{ccccccccccccccccccccccccccccccccccccccccccccccccccccc}
  &  & \  \downarrow (i+1)^{th} column \ \ \ \ \ \ &  &
  \end{array}
\]

$\nu^{j}_{\sigma^{1}[i] } := \left( \begin{array}{ccccccc}
     (n-j), & -j, & \ldots & ,-j, & (n-j), & \ldots & ,(n-j)  \\
   \end{array}
   \right),
$

\[\begin{array}{cccccccccccccccccccccc}
     &  &  &  &  &  &  &  &\downarrow (i+2)^{th} column &
  \end{array}
\]

$\nu^{j}_{\sigma^{2}[i] } := \left( \begin{array}{cccccccc}
     (n-j), & (n-j)& ,-j, & \ldots & ,-j, & (n-j), & \ldots & ,(n-j)  \\
   \end{array}
   \right),
$

\[\ldots \ldots \ldots \ldots \ldots \ldots \ldots \ldots \ldots 
\ldots \ldots \ldots \ldots \ldots \ldots \ldots \ldots
\ldots \ldots \ldots \ldots \ldots \ldots \ldots \ldots \ldots
\ldots \ldots \ldots \ldots\]

\[\begin{array}{ccccccccccccccccccccccc}
     \ \ \downarrow (i-2)^{th} column \ \ \ \ \downarrow (n-2)^{th} column \ \ \ \ \
  \end{array}
\]

$\nu^{j}_{\sigma^{n-2}[i] } := \left(
\begin{array}{ccccccccccccccc}
     -j, \ \ldots \ , -j , \ (n-j) , & \ldots & , (n-j) , \ -j , \ -j  \\
   \end{array}
   \right),
$

\[\begin{array}{ccccccccccccccccccccccccccccccc}
   \ \ \ \ \ \ \ \ \ \downarrow (i-1)^{th} column \ \ \ \ \downarrow 
   (n-1)^{th} column \ \ \ \ \ 
  \end{array}
\]

$\nu^{j}_{\sigma^{n-1}[i] } := \left( \begin{array}{cccccccc}
     -j, & \ldots & , -j , \ (n-j), & \ldots & ,(n-j) , \ -j  \\
   \end{array}
   \right),\\
$

where $\sigma^{k}[i]:=\{\sigma^{k}(1), \ldots , \sigma^{k}(i)\}$
for all $1\leq i \leq n-2,$ $1\leq j \leq n-1$ and $0\leq k \leq n-1.$\\

$\it{Remark}:$ \ It is easy to see $(\cite{SA})$ \ that for any 
$1\leq i \leq n-2$ \ and $0\leq t \leq n-1$   we have 
$\nu^{n-i-1}_{\sigma^{t}[i]} = \nu_{\sigma^{t}[i]}.$\\

If $0\neq a\in \mathbb{R}^{n},$ then $H_{a}$ will denote the
hyperplane of $\mathbb{R}^{n}$ through the origin with normal vector
a, that is,\[H_{a}=\{x\in \mathbb{R}^{n}\ | \ \langle x,a\rangle=0\},\] where
$\langle , \rangle$ is the usual inner product in $\mathbb{R}^{n}.$ The two closed
halfspaces bounded by $H_{a}$ are: \[H^{+}_{a}=\{x\in
\mathbb{R}^{n}\ | \ \langle x,a\rangle \geq 0\} \ and \  H^{-}_{a}=\{x\in
\mathbb{R}^{n}\ | \ \langle x,a\rangle \leq  0\}.\]

We will denote by $H_{\nu^{j}_{\sigma^{k}[i]}}$ the hyperplane of
$\mathbb{R}^{n}$ through the origin with normal vector
$\nu^{j}_{\sigma^{k}[i]}$, that is, \[H_{\nu^{j}_{\sigma^{k}[i]}}=\{x\in
\mathbb{R}^{n}\ | \ \langle x,\nu^{j}_{\sigma^{k}[i]}\rangle=0\},\] 
for all $1\leq i \leq n-2,$ $1\leq j \leq n-1$ and $0\leq k \leq n-1.$

Let $A\subset \mathbb{R}^{n},$ we denote
by \emph{aff}$(A)$ the affine space generated by $A.$ There is a unique linear 
subspace $V$ of $\mathbb{R}^{n}$ such that $\emph{aff}(A)=x_{0}+ V,$ for some
$x_{0}\in \mathbb{R}^{n}.$ The dimension of \emph{aff}$(A)$ is 
$dim(\emph{aff}(A))=dim_{\mathbb{R}}(V).$

Recall that a $polyhedral\ cone\ Q \subset \mathbb{R}^{n}$ is the
intersection of a finite number of closed subspaces of the form
$H^{+}_{a}.$ If $Q=H_{a_{1}}^{+}\cap \ldots \cap H_{a_{m}}^{+}$ is a polyhedral
cone, then \emph{aff}$(Q)$ is the intersection of those hyperplanes $H_{a_{i}}$ , 
$i= 1,\ldots,m$, that contain $Q$.(see {\cite[Proposition 1.2.] {BG}})
The $dimension \ of \ Q$ is the dimension of \emph{aff}$(Q)$, 
$dim(Q)=dim(\emph{aff}(Q))$.\\
If $A=\{\gamma_{1},\ldots, \ \gamma_{r}\}$ is a
finite set of points in $\mathbb{R}^{n}$ the $cone$ generated by
$A$, denoted by $\mathbb{R}_{+}{A},$  is defined as
\[\mathbb{R}_{+}{A}=\{\sum_{i=1}^{r}a_{i}\gamma_{i}\ | \ a_{i}\in
\mathbb{R}_{+}\ \mbox{for all} \ 1\leq i \leq n\}.\]
An important fact is that $Q$ is a polyhedral cone in
$\mathbb{R}^{n}$ if and only if there exists a finite set $A\subset
\mathbb{R}^{n}$ such that $Q=\mathbb{R}_{+}{A}$ (see {\cite{BG} or
\cite[Theorem 4.1.1.]{W}}).

Next we give some important definitions and results (see \cite{B},
 \cite{BH}, \cite{BG}, \cite{MS}, \cite{V}).

\begin{definition}
A \emph{proper face} of a polyhedral cone is a subset $F\subset \ Q$
 such
that there is a supporting hyperplane $H_{a}$ satisfying:

$1)$ $F=Q\cap H_{a}\neq \emptyset$,

$2)$ $Q\nsubseteq H_{a}$ and $Q\subset H^{+}_{a}$.
\end{definition}

The $dimension$ of a proper face $F$ of a polyhedral cone $Q$ is: 
$dim(F)=dim(\emph{aff}(F)).$

\begin{definition}
A cone $C$ is  \emph{pointed} if $0$ is a face of $C.$ \ Equivalently
we can require that $x\in C$ \ and $-x\in C$ \ $\Rightarrow$ \ $x=0.$
\end{definition}

\begin{definition}
The 1-dimensional faces of a pointed cone are called $extremal \
rays.$
\end{definition}

\begin{definition}
A proper face $F$ of a polyhedral cone $Q\subset \ \mathbb{R}^{n}$
is called a $facet$ of $Q$ if $\dim(F)=\dim(Q)-1.$
\end{definition}

\begin{definition}
If a polyhedral cone $Q$ is written as \[Q=H^{+}_{a_{1}}\cap \ldots
\cap H^{+}_{a_{r}}\] such that no  $H^{+}_{a_{i}}$ can be
omitted, then we say that this is an \emph{irreducible representation}
 of
$Q.$
\end{definition}

\begin{theorem}
Let $Q\subset \mathbb{R}^{n},$ \ $Q\neq \mathbb{R}^{n}$, be a
polyhedral cone with $dim(Q)=n.$ Then the halfspaces 
$H^{+}_{a_{1}}, \ldots , H^{+}_{a_{m}}$
in an irreducible representation $Q=H^{+}_{a_{1}}\cap \ldots 
\cap H^{+}_{a_{m}}$
are uniquely determined. In fact, the sets $F_{i}=Q 
\cap H_{a_{i}}, \  i = 1,\ldots, n,$ are the facets of
$Q$.
\end{theorem}
\begin{proof}
See {\cite[Theorem 1.6.] {BG}}
\end{proof}

\begin{definition}
Let $Q$ be a polyhedral cone in $\mathbb{R}^{n}$ \ with $\dim \ Q=n$
 and such that $Q\neq\mathbb{R}^{n}.$  Let 
\[Q=H^{+}_{a_{1}}\cap\ldots \cap H^{+}_{a_{r}}
\] 
be the irreducible representation of
$Q.$  If $a_{i}=(a_{i1},\ldots,a_{in}),$  then we call
\[H_{a_{i}}(x):=a_{i1}x_{1}+\ldots + a_{in}x_{n}=0,\ \ i\in [r],
\]
the \emph{equations of the cone} $Q.$
\end{definition}

\begin{definition}
The \emph{relative interior} $ri(Q)$ of a polyhedral cone is the 
interior of $Q$ with respect to the embedding of $Q$ into its affine
space \emph{aff}($Q$), in which $Q$ is full-dimensional.
\end{definition}

The following result gives us the description of the relative interior
of a polyhedral cone when we know its irreducible representation.

\begin{theorem}
Let $Q\subset \mathbb{R}^{n},$ \ $Q\neq \mathbb{R}^{n}$, be a
 polyhedral cone with $\dim(Q)=n$  and let 
\[(*) \ Q=H^{+}_{a_{1}}\cap \ldots \cap H^{+}_{a_{m}}
\]  
be an irreducible representation of $Q$  with $H^{+}_{a_{1}}, \ldots ,
 H^{+}_{a_{n}}$ pairwise distinct, 
where $a_{i} \in \mathbb{R}^{n} \setminus \{0\}$  for all $i.$  Set
 $F_{i}=Q \cap H_{a_{i}}$ \ for $i \in [r]$.  Then:\\
$a)$ \ $ri(Q)=\{x \in \mathbb{R}^{n} \ | \ \langle x, a_{1}\rangle \ >
 0, \ldots ,\langle x, a_{r}\rangle \ > 0 \}$,  
where $ri(Q)$ \ is the relative interior of $Q,$ \ which in this case
 is just the interior.\\
$b)$ \ Each facet $F$ \ of $Q$ \ is of the form $F=F_{i}$ \ for some
 $i.$\\
$c)$ \ Each $F_{i}$ \ is a facet of $Q$.
\end{theorem}
\begin{proof}
See  \cite[Theorems   8.2.15 and    3.2.1]{B}.
\end{proof}

\begin{theorem}{\bf(Danilov, Stanley)}
Let $R=K[x_{1},\ldots,x_{n}]$ be a polynomial ring over a field $K$ and
 $F$ a finite set
of monomials in $R.$
If $K[F]$ is normal, then the canonical module $\omega_{K[F]}$ of
 $K[F],$ with respect to
standard grading, can be expressed as an ideal of $K[F]$ generated by
 monomials
\[\omega_{K[F]}=(\{x^{a}| \ a\in {\mathbb{N}}A\cap
 ri({\mathbb{R_{+}}}A)\}),\] where $A=\log (F)$ and
$ri({\mathbb{R_{+}}}A)$ denotes the relative interior of
${\mathbb{R_{+}}}A.$
\end{theorem}

The formula above represents the canonical module of $K[F]$ \ as an
 ideal of $K[F]$ \ generated by monomials. For a comprehensive treatment of
 the Danilov-Stanley formula see \cite{BH},  \cite{MS} or \cite{V} . \\
 
\section{Polymatroids}
\[\]

Let $K$ be an infinite field, $n$ and $m$ be positive integers,
$[n]=\{1, 2, \ldots , n\}$. A nonempty finite set $B$ of ${\mathbb
N}^{n}$ is the \emph{base set of a discrete polymatroid}
${\bf{\mathcal{P}}}$ if for every $u=(u_{1}, u_{2}, \ldots,
u_{n})$, $v=(v_{1}, v_{2}, \ldots, v_{n})$ $\in B$ one has $u_{1}
+ u_{2} + \ldots + u_{n}=v_{1} + v_{2} + \ldots + v_{n}$ and for
all $i$ such that $u_{i}>v_{i}$ there exists $j$ such that
$u_{j}<v_{j}$ and $u+e_{j}-e_{i} \in B $, where $e_{k}$ denotes
the $k^{th}$ vector of the standard basis of ${\mathbb N}^{n}$. The
notion of discrete polymatroid is a generalization of the
classical notion of matroid, see \cite{E},  \cite{HH}, \cite{HHV}, 
\cite{O}, \cite{V1}, \cite{W1}. 
Associated with the base $B$ of a discrete polymatroid
${\bf{\mathcal{P}}}$ one has a $K-$algebra $K[ B ]$, called the
\emph{base ring} of ${\bf{\mathcal{P}}}$, defined to be the
$K-$subalgebra of the polynomial ring in $n$ indeterminates
$K[x_{1}, x_{2}, \ldots , x_{n}]$ generated by the monomials
$x^{u}$ with $u \in B$. From \cite{HH} the algebra $K[ B ]$ is
known to be normal and hence Cohen-Macaulay.

If $A_{i}$ are some nonempty subsets of $[n]$ for $1\leq i\leq
m$, ${\bf{\mathcal{A}}}=\{A_{1},\ldots,A_{m}\}$, then the
set of the vectors $\sum_{k=1}^{m} e_{i_{k}}$ with $i_{k} \in
A_{k}$ is the base of a polymatroid, called the \emph{transversal
polymatroid presented} by ${\bf{\mathcal{A}}}.$ The \emph{base ring} of
 a
transversal polymatroid presented by ${\bf{\mathcal{A}}}$ 
is the ring 
\[K[{\bf{\mathcal{A}}}]:=K[x_{i_{1}}\cdot\cdot\cdot
x_{i_{m}}\ | \  i_{j}\in A_{j},1\leq j\leq m].\]

\section{The type of base ring associated to transversal polymatroids 
with the cone of dimension $n$ with $n+1$ facets.}

\begin{lemma}
Let $n\in \mathbb{N},$ \ $n \geq 3,$ \ $1\leq i\leq n-2$ \ 
and $1\leq j\leq n-1.$ \ We consider the following two cases:

$a)$ \ If $i+j\leq n-1,$ \ then let $A:=\{
log(x_{j_{1}}\cdot\cdot\cdot x_{j_{n}}) \ | \
j_{k}\in A_{k},\ for \ all \ 1\leq k \leq n
\}\subset \mathbb{N}^{n}$ the exponent set of generators of
$K-$algebra $K[{\bf{\mathcal{A}}}],$ where
${\bf{\mathcal{A}}}=\{A_{1}=[n],\ldots,A_{i}=[n],
A_{i+1}=[n]\setminus [i],\ldots,A_{i+j}=[n]\setminus
[i],A_{i+j+1}=[n],\ldots,A_{n}=[n]\}$.

$b)$ \ If $i+j\geq n,$ \ then let $A:=\{
log(x_{j_{1}}\cdot\cdot\cdot x_{j_{n}}) \ | \
j_{k}\in A_{k},\ for \ all \ 1\leq k \leq n
\}\subset \mathbb{N}^{n}$ the exponent set of generators of
$K-$algebra $K[{\bf{\mathcal{A}}}],$ where
${\bf{\mathcal{A}}}=\{A_{1}=[n]\setminus
[i],\ldots,A_{i+j-n}=[n]\setminus
[i],
A_{i+j-n+1}=[n],\ldots,A_{i}=[n], A_{i+1}=[n]\setminus
[i],\ldots,A_{n}=[n]\setminus
[i]\}$.\\
Then the cone generated by $A$ has the
irreducible representation:
\[\mathbb{R}_{+}A= \bigcap_{a\in N}H^{+}_{a},\] where
$N=\{\nu^{j}_{\sigma^{0}[i]}, \ e_{k} \ | \ 1\leq k \leq
n\},$ \  $\{e_{i}\}_{1 \leq i \leq n}$ \ is the
canonical base of $\mathbb{R}^{n}.$
\end{lemma}

\begin{proof}
We denote by 
$
J_{k}=\left \{ \begin{array}[pos]{cccccccccccccccc}
(n-j) \ e_{k}+ j \ e_{i+1} \ , \ \ if \ \ 1\leq k \leq i \ \ \ \ \ \ \\

(n-j) \ e_{1}+j \ e_{r} \ , \ \ \ \ \ if \ i+2\leq r \leq n \ \
\end{array} \right .
$ and by $J=n \ e_{n}.$ \
Since $A_{t}=[n]$ \ for any $t \in \{1,\ldots, i\}\cup \{i+j+1,
\ldots, n\}$ \ and $A_{r}=[n] \setminus [i]$ \ for any $r \in 
\{i+1,\ldots, i+j\}$ it is easy to see that for any $k \in 
\{1,\ldots, i\}$ \ and  $r \in \{i+2,\ldots, n\}$ \ the set of 
monomials $x^{n-j}_{k} \ x^{j}_{i+1},$ \ $x^{n-j}_{1} \ x^{j}_{r},$ 
\ $x^{n}_{n}$ \ are a subset of the generators of $K-$algebra 
$K[{\bf{\mathcal{A}}}].$ \ Thus the set \[\{J_{1}, \ldots, J_{i},
J_{i+2}, \ldots, J_{n}, J\}\subset A.\]
If we denote by $C$ \ the matrix with the rows the coordinates of 
the $\{J_{1}, \ldots, J_{i}, J_{i+2}, \ldots, J_{n}, J\},$ \ 
then by a simple computation we get $|det\ (C)|=n \ (n-j)^{i} \ 
j^{n-i-1}$ \ for any $1\leq i\leq n-2$ \ and $1\leq j\leq n-1.$ \ 
Thus, we get that the set \[\{J_{1}, \ldots, J_{i}, J_{i+2}, 
\ldots, J_{n}, J\}\subset A.\] is linearly independent and it 
follows that $dim \ \mathbb{R}_{+}A = n.$ \ Since $\{J_{1}, 
\ldots, J_{i}, J_{i+2}, \ldots, J_{n}\}$ is linearly independent 
and lie on the hyperplane $H_{\nu^{j}_{\sigma^{0}[i]}}$ we have 
that $dim(H_{\nu^{j}_{\sigma^{0}[i]}}\cap \mathbb{R}_{+}A)=n-1.$\\
Now we will prove that $\mathbb{R}_{+}A \subset H_{a}^{+}$ \ for 
all $a \in N.$ \ It is enough to show that for all vectors 
$P \in A,$ \ $\langle P,a \rangle \geq 0$ for all $a \in N.$ \ 
Since $\{e_{k}\}_{1 \leq k \leq n}$ \ is  the
canonical base of $\mathbb{R}^{n},$\ we get that 
$\langle P,\ e_{k}\rangle \geq 0$ for any $1\leq k \leq n.$  \ 
Let $P \in A,$ \ $P=log(x_{j_{1}}\cdot\cdot\cdot x_{j_{n}}).$ \ 
We have two possibilities:\\ 
$a)$ \ If $i+j\leq n-1,$ \ then let $t$ 
 be the number of $j_{k},$ \ such that \ $k \in \{1,\ldots, 
i\}\cup\{i+j+1,\ldots, n\}$ \ and $j_{k} \in [i].$ \ Thus $t \leq n-j.$ \ 
Then $\langle P, \nu^{j}_{\sigma^{0}[i]}\rangle =-t \ j + (n-j-t)\ 
(n-j) + j(n-j)= n(n-j-t) \geq 0.$\\  
$b)$ \ If $i+j\geq n,$ \ 
then let $t$ be the number of $j_{k},$ \ such that \ $i+j-n+1
\leq k \leq i$ \ and $j_{k} \in [i].$ \ Thus $t \leq n-j.$ \ 
Then $\langle P, \nu^{j}_{\sigma^{0}[i]}\rangle=-t \ j + (n-j-t)\ (n-j) 
+ j(n-j)= n(n-j-t) \geq 0.$\\ Thus \[\mathbb{R}_{+}A\subseteq 
\bigcap_{a\in N}H^{+}_{a}.\]Now we will prove the converse inclusion: 
$\mathbb{R}_{+}A\supseteq \bigcap_{a\in N}H^{+}_{a}.$\\ It is  
enough to prove that the extremal rays of the cone 
$\bigcap_{a\in N}H^{+}_{a}$ \ are in
${\mathbb{R}_{+}}A.$ \ Any extremal ray of the cone  
$\bigcap_{a\in N}H^{+}_{a}$ \ can be written  as the intersection of 
$n-1$ hyperplanes $H_{a}$ \ , with $a\in N.$ \ 
There are two possibilities to obtain extremal rays by intersection of 
$n-1$ hyperplanes.\\ 
$First \ case.$\\ Let \ $1\leq i_{1} < \ldots < 
i_{n-1}\leq n$ \ be a sequence and $\{t\}=[n]\setminus \{i_{1},
\ldots,i_{n-1}\}.$ \ The system of equations:
$ \ (*) \ \begin{cases}
z_{i_{1}}=0,\\
\vdots \\
z_{i_{n-1}}=0
\end{cases}$ admits the solution $x \in \mathbb{Z}_{+}^{n},$ \ $
x=\left( \begin{array}[pos]{c}
    x_{1}\\
    \vdots\\
    x_{n}
\end{array} \right)
$ with $ | \ x \ |=n,$ \ $x_{k}=n\cdot\delta_{kt}$ \ for all 
$1\leq k \leq n,$ where $\delta_{kt}$ is Kronecker symbol.\\
There are two possibilities:\\ $1)$ \ If $1\leq t \leq i,$ \ 
then $H_{\nu^{j}_{\sigma^{0}[i]}}(x) < 0$ \ and thus $x\notin 
\bigcap_{a\in N}H^{+}_{a}.$\\
$2)$ \ If $i+1\leq t \leq n,$ \ then $H_{\nu^{j}_{\sigma^{0}[i]}}(x) > 0$ \ 
and thus $x\in \bigcap_{a\in N}H^{+}_{a}$ and  is an extremal ray.\\ 
Thus, there exist $n-i$ \ sequences \ $1\leq i_{1}< \ldots < i_{n-1}\leq n$ \ 
such that the system of equations $(*)$ \ has a solution 
$x \in \mathbb{Z}_{+}^{n}$ \ with $ | \ x \ |=n$ \ 
and $H_{\nu^{j}_{\sigma^{0}[i]}}(x) > 0.$\\ The extremal rays are: 
$\{n \ e_{k} \ | \ i+1 \leq k \leq n\}.$\\
$Second \ case.$\\ 
Let \ $1\leq i_{1}< \ldots < i_{n-2}\leq n$ \ be a sequence and 
$\{r, \ s\}=[n]\setminus \{i_{1},\ldots,i_{n-2}\},$ \ with $r < s.$\\
Let $x \in \mathbb{Z}_{+}^{n},$ \ with $ | \ x \ |=n,$ \ 
a solution of system of equations:\\
$ \ (**) \ \begin{cases}
z_{i_{1}}=0,\\
\vdots \\
z_{i_{n-2}}=0,\\
-j \ z_{1}-\ldots-j \ z_{i}+(n-j)z_{i+1}+\ldots+(n-j)z_{n}=0.
\end{cases}$\\
There are two possibilities:\\
$1)$ \ If $1\leq r \leq i$ and $i+1\leq s
 \leq n,$ then the system of equations $(**)$ admits the solution
$
x=\left( \begin{array}[pos]{c}
    x_{1}\\
    \vdots\\
    x_{n}
\end{array} \right) \in \mathbb{Z}_{+}^{n},$ \ with $ | \ x \ |=n,$ \ 
with $x_{t}=j \ \delta_{ts}+(n-j)\delta_{tr}$ \ for all $1\leq t \leq n.$\\
$2)$ If $1\leq r,s \leq i$ \ or $i+1\leq r,s \leq n,$ \ then there exist 
no solution $x\in \mathbb{Z}_{+}^{n}$ \ for the system of equations $(**)$ 
\ because otherwise $H_{\nu^{j}_{\sigma^{0}[i]}}(x) > 0$ \ or 
$H_{\nu^{j}_{\sigma^{0}[i]}}(x) < 0.$\\
Thus, there exist $i(n-i)$ \ sequences \ $1\leq i_{1}< \ldots < i_{n-2}\leq n$ 
such that the system of equations $(**)$ \ has a solution 
$x \in \mathbb{Z}_{+}^{n}$ \ with $ | \ x \ |=n$ \ and the extremal rays are: 
\[\{(n-j)e_{r}+j \ e_{s} \ | \ 1 \leq r \leq i \ and \ i+1 \leq s \leq n\}.\]\\
In conclusion, there exist $(i+1)(n-i)$ \ extremal rays of the cone 
${\bigcap}_{a\in N}H_{a}^{+}$:\[R:=\{n e_{k} \ | \ i+1 \leq k \leq n\}
\cup\{(n-j)e_{r}+j \ e_{s} \ | \ 1 \leq r \leq i \ and \ i+1 \leq s \leq n\}.\]
Since  $A=\{log(x_{j_{1}}\cdot\cdot\cdot x_{j_{n}}) \ | \
j_{k}\in A_{k} \in {\bf{\mathcal{A}}},\ for \ all \ 1\leq k \leq n
\}\subset \mathbb{N}^{n}$ \ where ${\bf{\mathcal{A}}}=\{A_{1}=[n],\ldots,A_{i}=[n],
A_{i+1}=[n]\setminus [i],\ldots,A_{i+j}=[n]\setminus
[i],A_{i+j+1}=[n],\ldots,A_{n}=[n]\}$ \ or ${\bf{\mathcal{A}}}=\{A_{1}=[n]\setminus
[i],\ldots,A_{i+j-n}=[n]\setminus [i],
A_{i+j-n+1}=[n],\ldots,A_{i}=[n], A_{i+1}=[n]\setminus
[i],\ldots,A_{n}=[n]\setminus
[i]\},$ \ it is easy to see that $R\subset A.$ \ Thus, we have $\mathbb{R}_{+}A= 
\bigcap_{a\in N}H^{+}_{a}.$\\ The representation is
irreducible because if we delete, for some $k,$ \ the hyperplane with
the normal $e_{k}$, then the $k^{'th}$ coordinate
of $log(x_{j_{1}}\cdot\cdot\cdot x_{j_{n}})$
would be negative, which it is impossible; and if we delete the
hyperplane with the normal $\nu^{j}_{\sigma^{0}[i]},$ then the cone
${\mathbb{R_{+}}}A$ would be presented by $A=\{
log(x_{j_{1}}\cdot\cdot\cdot x_{j_{n}}) \ | \ j_{k}\in [n], \ for \ all \ 1\leq
k \leq n \}$ which it is impossible. \ Thus the representation
${\mathbb{R_{+}}}A= \bigcap_{a\in N}H^{+}_{a}$ is irreducible.
\end{proof}

\begin{lemma}
Let $n\in \mathbb{N},$ \ $n \geq 3,$ \ $t\geq 1,$ \ $1\leq i\leq n-2$ \ and 
$1\leq j\leq n-1.$ \ Let $A=\{ log(x_{j_{1}}\cdot\cdot\cdot x_{j_{n}}) \ | \
j_{\sigma^{t}{(k)}}\in A_{\sigma^{t}{(k)}}, 1\leq k \leq
n\}\subset \mathbb{N}^{n}$ the exponent set of generators of
$K-$algebra  $K[{\bf{\mathcal{A}}}],$ where ${\bf{\mathcal{A}}}=\{
A_{\sigma^{t}{(k)}} \ | \ A_{\sigma^{t}{(k)}}=[n],\ for \ k \in [i] \cup 
([n]\setminus [i+j]) \ and \ A_{\sigma^{t}{(k)}}=[n]\setminus \sigma^{t}[i], 
\ for \ i+1 \leq k \leq i+j\},$ \  if $i+j \leq n-1$ or ${\bf{\mathcal{A}}}=\{
A_{\sigma^{t}{(k)}} \ | \ A_{\sigma^{t}{(k)}}=[n]\setminus [i],\ for \ 
k \in [i+j-n] \cup ([n]\setminus [i]) \ and \ A_{\sigma^{t}{(k)}}=[n], \ 
for \ i+j-n+1 \leq k \leq i\},$ \ if $i+j \geq n.$  Then the cone generated by $A$
has the irreducible representation:
\[{\mathbb{R_{+}}}A= \bigcap_{a\in N}H^{+}_{a},\] where 
$N=\{\nu^{j}_{\sigma^{t}[i]},\ e_{k} \ | \ 1\leq k \leq
n\}.$
\end {lemma}
For the proof note that the algebras from Lemma 4.1. and 4.2.
are isomorphic.

\begin{lemma}
$a)$ Let $1\leq t \leq n-i-j-1$, $s\geq 2$ and 
$\beta \in {\bf{\mathbb{N}}}^{n}$ such that 
$H_{\nu^{j}_{\sigma^{0}[i]}}(\beta)= n(n-i-j-t)$ 
and $| \ \beta \ |=sn$. 
Then $\beta_{1}+\ldots+\beta_{i}=(n-j)(s-1)+i+t$ and
$\beta_{i+1}+\ldots+\beta_{n}=n+j(s-1)-i-t$.\\
$b)$ Let $r\geq 2$, $s\geq r$, $1\leq t\leq r(n-j)-i$ and 
$\beta \in {\bf{\mathbb{N}}}^{n}$ such that 
$H_{\nu^{j}_{\sigma^{0}[i]}}(\beta)= nt$ and $| \ \beta \ |=sn$. 
Then $\beta_{1}+\ldots+\beta_{i}=(n-j)s-t$ and
$\beta_{i+1}+\ldots+\beta_{n}=js+t$.
\end{lemma}
\begin{proof}
$a)$ Let $c=\beta_{1}+\ldots+\beta_{i}$ and thus $\beta_{i+1}+\ldots+\beta_{n}=sn-c.$
Then $H_{\nu^{j}_{\sigma^{0}[i]}}(\beta)=-jc+(n-j)(sn-c)=n(sn-c-js)$. Since
$H_{\nu^{j}_{\sigma^{0}[i]}}(\beta)=n(n-i-j-t)$ it follows that
$sn-c-js=n-i-j-t$ and so $c=(n-j)(s-1)+i+t$. 
Thus, $\beta_{1}+\ldots+\beta_{i}=(n-j)(s-1)+i+t$ and
$\beta_{i+1}+\ldots+\beta_{n}=sn-c=n+j(s-1)-i-t$.\\
$b)$ The proof goes as that for $a).$
\end{proof}

The main result of this paper is the following:

\begin{theorem}
Let $R=K[x_{1},\ldots,x_{n}]$ be a polynomial ring standard graded  
over a field $K$ \ and $A$ the finite set
of monomials in $R$ such that the hypothesis of $Lemma \ 4.1.$  
are satisfied.  Then: \\
$a)$ \ If $i+j\leq n-1,$  then the type of $K[{\bf{\mathcal{A}}}]$ is:
\[type(K[{\bf{\mathcal{A}}}])=1+ \sum_{t=1}^{n-i-j-1} \binom{n+i-j+t-1}
{i-1}\binom{n-i+j-t-1}{n-i-1}.\]
$b)$  \ If $i+j\geq n,$ \ then the type of $K[{\bf{\mathcal{A}}}]$ is:
\[type(K[{\bf{\mathcal{A}}}])=\sum^{r(n-j) -i}_{t=1}
\binom{r(n-j)-t-1}{i-1}\binom{rj+t-1}{n-i-1},\]
where $r=\left\lceil \frac{i+1}{n-j} \right\rceil$ $(\left\lceil x \right\rceil 
is \ the \ least \ integer \ \geq x).$
\end{theorem}

\begin{proof}
Since $K[{\bf{\mathcal{A}}}]$ is normal (\cite{HH}), 
then the canonical module $\omega_{K[{\bf{\mathcal{A}}}]}$ of 
$K[{\bf{\mathcal{A}}}],$ 
with respect to standard grading, can be expressed as an ideal of 
$K[{\bf{\mathcal{A}}}]$ generated by monomials
\[\omega_{K[{\bf{\mathcal{A}}}]}=(\{x^{a}| \ a\in {\mathbb{N}}A\cap 
ri({\mathbb{R_{+}}}A)\}),\] 
where $A$ is the exponent set of $K-$ algebra $K[{\bf{\mathcal{A}}}]$
and $ri({\mathbb{R_{+}}}A)$ denotes the relative interior of
${\mathbb{R_{+}}}A.$ By Lemma 4.1.,  the cone generated by $A$ has the
irreducible representation:
\[\mathbb{R}_{+}A= \bigcap_{a\in N}H^{+}_{a},\] where
$N=\{\nu^{j}_{\sigma^{0}[i]}, \ e_{k} \ | \ 1\leq k \leq
n\},$ \  $\{e_{k}\}_{1 \leq k \leq n}$ \ being the
canonical base of $\mathbb{R}^{n}.$\\
$a)$ \ Let $i\in [n-2],$ $j\in [n-1]$ 
such that $i+j\leq n-1,$ \ $u_{t}=n+i-j+t,$ $v_{t}=n-i+j-t$ for 
any $t\in [n-i-j-1].$  We will denote by $M$ the set:
\[M=\{\alpha \in \mathbb{N}^{n}\ | \ \alpha_{k}\geq 1, \ |(\alpha_{1},\ldots,
\alpha_{i})|=u_{t}, \ |(\alpha_{i+1},\ldots,\alpha_{n})|=v_{t} 
\ for \ any \ k\in [n], \ i\in [n-2], \ \]\[ j\in [n-1] \  
 such \  that \ i+j\leq n-1 \ and \ t\in [n-i-j-1]\}.\]
We will show that : \[{\mathbb{N}}A\cap ri({\mathbb{R_{+}}}A)=
((1,\ldots,1)+({\mathbb{N}}A\cap{\mathbb{R_{+}}}A))\cup
\bigcup_{\alpha \in M} ( \{\alpha\}+({\mathbb{N}}A\cap{\mathbb{R_{+}}}A)).\]
Since for any $\alpha \in M,$ \ $H_{\nu^{j}_{\sigma^{0}[i]}}(\alpha)=
-j(n+i-j+t)+(n-j)(n-i+j-t)=n(n-i-j-t) \ > \ 0$ and $H_{\nu^{j}_
{\sigma^{0}[i]}}((1,\ldots,1))=n(n-i-j) > 0,$ \ it follows that : 
\[{\mathbb{N}}A\cap ri({\mathbb{R_{+}}}A) \supseteq
((1,\ldots,1)+({\mathbb{N}}A\cap{\mathbb{R_{+}}}A))\cup
\bigcup_{\alpha \in M} ( \{\alpha\}+({\mathbb{N}}A\cap{\mathbb{R_{+}}}A)).\]
Let $\beta \in {\mathbb{N}}A\cap ri({\mathbb{R_{+}}}A),$ then $\beta_{k}\geq 1$
for any $k\in [n].$ Since $H_{\nu^{j}_{\sigma^{0}[i]}}((1,\ldots,1))=n(n-i-j) > 0$
it follows that $(1,\ldots,1)\in ri({\mathbb{R_{+}}}A).$ Let $\gamma 
\in \mathbb{N}^{n}, \ \gamma=\beta - (1,\ldots,1).$ It is clear that
$H_{\nu^{j}_{\sigma^{0}[i]}}(\gamma)=H_{\nu^{j}_{\sigma^{0}[i]}}
(\beta)-n(n-i-j).$
If $H_{\nu^{j}_{\sigma^{0}[i]}}(\beta)\geq n(n-i-j),$ then  
$H_{\nu^{j}_{\sigma^{0}[i]}}(\gamma)\geq 0.$ Thus $\gamma \in 
{\mathbb{N}}A\cap {\mathbb{R_{+}}}A.$ \\
If $H_{\nu^{j}_{\sigma^{0}[i]}}(\beta)< n(n-i-j),$ \ 
then let $1\leq t \leq n-i-j-1$ such that 
$H_{\nu^{j}_{\sigma^{0}[i]}}(\beta)= n(n-i-j-t).$

We claim that for any $\beta \in {\mathbb{N}}A\cap ri({\mathbb{R_{+}}}A)$ 
with $| \ \beta \ |=sn\geq 2n$
and $t\in [n-i-j-1]$ such that  
$H_{\nu^{j}_{\sigma^{0}[i]}}(\beta)= n(n-i-j-t)$ we can find 
$\alpha \in M$ with $H_{\nu^{j}_{\sigma^{0}[i]}}(\alpha)= n(n-i-j-t)$
and $\beta-\alpha \in {\mathbb{N}}A\cap{\mathbb{R_{+}}}A.$\\
We proceed by induction on $s\geq 2.$
If $s=2$, then it is clear that $\beta \in M$. 
Indeed, then $\beta_{1}+\ldots+\beta_{i}=n-a,$ 
$\beta_{i+1}+\ldots+\beta_{n}=n+a$ for some 
$a\in \mathbb{Z}$ and so
$H_{\nu^{j}_{\sigma^{0}[i]}}(\beta)=n(n-2j+a)=
n(n-i-j-t)$ that is
$a=j-i-t$. It follows $\beta \in M.$\\
Suppose $s>2.$ Since $\beta_{1}+\ldots+\beta_{i}= 
(n-j)(s-1)+i+t\geq 2(n-j)+i+t$ by $Lemma \ 4.3$ \ 
we can get $\gamma_{e}$, with $0\leq \gamma_{e}\leq \beta_{e}$ 
for any $1\leq e \leq i$,
such that $| \ (\gamma_{1},\ldots,\gamma_{i}) \ |=n-j.$
Since $\beta_{i+1}+\ldots+\beta_{n}=j(s-1)+(n-i-t)\geq 2j+n-i-t$ 
by $Lemma \ 4.3$ \ 
we can get $\gamma_{e}$, with $0\leq \gamma_{e}\leq \beta_{e}$ 
for any $i+1\leq e \leq n$,
such that $| \ (\gamma_{i+1},\ldots,\gamma_{n}) \ |=j.$
It is clear that $\beta^{'}=\beta-\gamma \in \mathbb{R}_{+}^{n}$, 
\[ | \ (\beta_{1}^{'},\ldots,\beta_{i}^{'}) \ |=
\sum_{e=1}^{i}(\beta_{e}-\gamma_{e})=(n-j)(s-2)+i+t,\]
\[ | \ (\beta_{i+1}^{'},\ldots,\beta_{n}^{'})  \ |=
\sum_{e=i+1}^{n}(\beta_{e}-\gamma_{e})=n+j(s-2)-i-t\] and
\[H_{\nu^{j}_{\sigma^{0}[i]}}(\beta^{'})=
H_{\nu^{j}_{\sigma^{0}[i]}}(\beta)-H_{\nu^{j}_{\sigma^{0}[i]}}(\gamma)=
n(n-i-j-t)-[(-j)(n-j)+(n-j)j]=n(n-i-j-t).\]
So, $\beta^{'}\in {\mathbb{N}}A\cap ri({\mathbb{R_{+}}}A)$ and since 
$| \ \beta^{'} \ |=n(s-1)$ we get by induction hypothesis that 
\[\beta^{'}\in \bigcup_{\alpha \in M} ( \{\alpha\}+({\mathbb{N}}A\cap{\mathbb{R_{+}}}A)).\]
Since $H_{\nu^{j}_{\sigma^{0}[i]}}(\gamma)=0$ and $\gamma \in \mathbb{R}_{+}^{n}$
we get that $ \gamma \in {\mathbb{N}}A \cap {\mathbb{R_{+}}}A$ and so
\[\beta \in \bigcup_{\alpha \in M} ( \{\alpha\}+({\mathbb{N}}A\cap{\mathbb{R_{+}}}A)).\]

Thus:
\[{\mathbb{N}}A\cap ri({\mathbb{R_{+}}}A) \subseteq
((1,\ldots,1)+({\mathbb{N}}A\cap{\mathbb{R_{+}}}A))\cup
\bigcup_{\alpha \in M} ( \{\alpha\}+({\mathbb{N}}A\cap{\mathbb{R_{+}}}A)).\]
So, the canonical module $\omega_{K[{\bf{\mathcal{A}}}]}$ of 
$K[{\bf{\mathcal{A}}}],$ 
with respect to standard grading, can be expressed as an ideal of 
$K[{\bf{\mathcal{A}}}],$ generated by monomials
\[\omega_{K[{\bf{\mathcal{A}}}]}=(\{x_{1}\cdot\cdot\cdot x_{n}, 
\ x^{\alpha}| \ \alpha\in M\}) K[{\bf{\mathcal{A}}}].\]
The type of $K[{\bf{\mathcal{A}}}]$ is the number of generator of
the canonical module. So, \ $type(K[{\bf{\mathcal{A}}}])=\#(M)+1,$ 
where $\#(M)$ is the cardinal of $M.$\\
Note that for 
 any $t\in [n-i-j-1],$ the equation: $\alpha_{1}+ \ldots +
\alpha_{i}=n+i-j+t$ has $\binom{n+i-j+t-1}{i-1}$ distinct nonnegative
integer solutions with $\alpha_{k}\geq 1,$ for any $k\in [i],$ respectively
$\alpha_{i+1}+ \ldots +
\alpha_{n}=n-i+j-t$ has $\binom{n-i+j-t-1}{n-i-1}$ distinct nonnegative
integer solutions with $\alpha_{k}\geq 1$ for any $k\in [n]\setminus[i].$
Thus, the cardinal of $M$ is \[\#(M)=\sum_{t=1}^{n-i-j-1} \binom{n+i-j+t-1}
{i-1}\binom{n-i+j-t-1}{n-i-1}.\] 

$b)$ \ Let $i\in [n-2],$ $j\in [n-1]$ 
such that $i+j\geq n,$ \ $u_{t}=r(n-j)-t,$ $v_{t}=rj+t$ for 
any $t\in [r(n-j)-i],$ where $r=\left\lceil \frac{i+1}{n-j} \right\rceil.$  
We will denote by $M^{'}$ the set:
\[M^{'}=\{\alpha \in \mathbb{N}^{n}\ | \ \alpha_{k}\geq 1, \ |(\alpha_{1},\ldots,
\alpha_{i})|=u_{t}, \ |(\alpha_{i+1},\ldots,\alpha_{n})|=v_{t} 
\ for \ any \ k\in [n], \ i\in [n-2], \ \]\[ j\in [n-1] \  
 such \  that \ i+j\geq n \ and \ t\in [r(n-j)-i]\}.\]
We will show that : \[{\mathbb{N}}A\cap ri({\mathbb{R_{+}}}A)=
\bigcup_{\alpha \in M^{'}} ( \{\alpha\}+({\mathbb{N}}A\cap{\mathbb{R_{+}}}A)).\]
Since for any $\alpha \in M^{'},$ \ $H_{\nu^{j}_{\sigma^{0}[i]}}(\alpha)=
-j(r(n-j)-t)+(n-j)(rj+t)=nt \ > \ 0,$ \ it follows that : 
\[{\mathbb{N}}A\cap ri({\mathbb{R_{+}}}A) \supseteq
\bigcup_{\alpha \in M^{'}} ( \{\alpha\}+({\mathbb{N}}A\cap{\mathbb{R_{+}}}A)).\]
Let $\beta \in {\mathbb{N}}A\cap ri({\mathbb{R_{+}}}A),$ then $\beta_{k}\geq 1$
for any $k\in [n].$ Since $H_{\nu^{j}_{\sigma^{0}[i]}}((1,\ldots,1))=n(n-i-j)<0$
it follows that $(1,\ldots,1)\notin ri({\mathbb{R_{+}}}A).$ We claim that
$| \ \beta \ |\geq rn.$ Indeed, since $\beta \in {\mathbb{N}}A\cap ri({\mathbb
{R_{+}}}A)$ and $| \ \beta \ | =sn,$ it follows that:
\[H_{\nu^{j}_{\sigma^{0}[i]}}(\beta)=-j\sum_{k=1}^{i}\beta_{i}+(n-j)
(sn-\sum_{k=1}^{i}\beta_{i}) > 0\Longleftrightarrow \sum_{k=1}^{i}\beta_{i}
< (n-j)s.\] Hence $i+1\leq s(n-j)$ and so
$r=\left\lceil \frac{i+1}{n-j} \right\rceil \leq s.$\\

We claim that for any $\beta \in {\mathbb{N}}A\cap ri({\mathbb{R_{+}}}A)$
with $| \ \beta \ |=sn \geq rn$ and $t\in [r(n-j)-i]$ such that 
$H_{\nu^{j}_{\sigma^{0}[i]}}(\beta)= nt$ we can find 
$\alpha \in M^{'}$ with $H_{\nu^{j}_{\sigma^{0}[i]}}(\alpha)= nt$
such that 
$\beta-\alpha \in {\mathbb{N}}A\cap{\mathbb{R_{+}}}A.$ \\
We proceed by induction on $s\geq r.$
If $s=r$, then it is clear that $\beta \in M^{'}$. 
Indeed, then $\beta_{1}+\ldots+\beta_{i}=n(r-1)-a,$ 
$\beta_{i+1}+\ldots+\beta_{n}=n+a$ for some 
$a\in \mathbb{Z}$ and so
$H_{\nu^{j}_{\sigma^{0}[i]}}(\beta)=n(n-jr+a)=
nt$ that is
$a=jr-n+t$. It follows $\beta \in M^{'}.$\\
Suppose $s>r$. Since $\beta_{1}+\ldots+\beta_{i}= 
(n-j)s-t\geq (n-j)(s-r)+i$ by $Lemma \ 4.3$
we can get $\gamma_{e}$, with $0\leq \gamma_{e}\leq \beta_{e}$ 
for any $1\leq e \leq i$,
such that $| \ (\gamma_{1},\ldots,\gamma_{i}) \ |=n-j.$
Since $\beta_{i+1}+\ldots+\beta_{n}=js+t\geq jr+t$ by $Lemma \ 4.3$
we can get $\gamma_{e}$, with $0\leq \gamma_{e}\leq \beta_{e}$ 
for any $i+1\leq e \leq n$,
such that $| \ (\gamma_{i+1},\ldots,\gamma_{n}) \ |=j.$
It is clear that $\beta^{'}=\beta-\gamma \in \mathbb{R}_{+}^{n}$, 
\[ | \ (\beta_{1}^{'},\ldots,\beta_{i}^{'}) \ |=
\sum_{e=1}^{i}(\beta_{e}-\gamma_{e})=(n-j)(s-1)-t,\]
\[ | \ (\beta_{i+1}^{'},\ldots,\beta_{n}^{'})  \ |=
\sum_{e=i+1}^{n}(\beta_{e}-\gamma_{e})=j(s-1)+t\] and
\[H_{\nu^{j}_{\sigma^{0}[i]}}(\beta^{'})=
H_{\nu^{j}_{\sigma^{0}[i]}}(\beta)-H_{\nu^{j}_{\sigma^{0}[i]}}(\gamma)=
nt-[(-j)(n-j)+(n-j)j]=nt.\]
So, $\beta^{'}\in {\mathbb{N}}A\cap ri({\mathbb{R_{+}}}A)$ and since 
$| \ \beta^{'} \ |=n(s-1)$ we get by induction hypothesis that 
\[\beta^{'}\in \bigcup_{\alpha \in M^{'}} ( \{\alpha\}+({\mathbb{N}}A\cap{\mathbb{R_{+}}}A)).\]
Since $H_{\nu^{j}_{\sigma^{0}[i]}}(\gamma)=0$ and $\gamma \in \mathbb{R}_{+}^{n}$
we get that $ \gamma \in {\mathbb{N}}A \cap {\mathbb{R_{+}}}A$ and so
\[\beta \in \bigcup_{\alpha \in M^{'}} ( \{\alpha\}+({\mathbb{N}}A\cap{\mathbb{R_{+}}}A)).\]

Thus:
\[{\mathbb{N}}A\cap ri({\mathbb{R_{+}}}A) \subseteq
\bigcup_{\alpha \in M^{'}} ( \{\alpha\}+({\mathbb{N}}A\cap{\mathbb{R_{+}}}A)).\]
So, the canonical module $\omega_{K[{\bf{\mathcal{A}}}]}$ of 
$K[{\bf{\mathcal{A}}}],$ 
with respect to standard grading, can be expressed as an ideal of 
$K[{\bf{\mathcal{A}}}],$ generated by monomials
\[\omega_{K[{\bf{\mathcal{A}}}]}=(\{x^{\alpha}| \ \alpha\in M^{'}\})
 K[{\bf{\mathcal{A}}}].\]
The type of $K[{\bf{\mathcal{A}}}]$ is the number of generator of
the canonical module. So, \ $type(K[{\bf{\mathcal{A}}}])=\#(M^{'}),$ 
where $\#(M^{'})$ is the cardinal of $M^{'}.$\\
Note that for 
 any $t\in [r(n-j)-i],$ the equation: $\alpha_{1}+ \ldots +
\alpha_{i}=r(n-j)-t$ has $\binom{r(n-j)-t-1}{i-1}$ distinct nonnegative
integer solutions with $\alpha_{k}\geq 1,$ for any $k\in [i],$ respectively
$\alpha_{i+1}+ \ldots +
\alpha_{n}=rj+t$ has $\binom{rj+t-1}{n-i-1}$ distinct nonnegative
integer solutions with $\alpha_{k}\geq 1$ for any $k\in [n]\setminus[i].$\\
Thus, the cardinal of $M^{'}$ is \[\#(M^{'})=\sum^{r(n-j) -i}_{t=1}
\binom{r(n-j)-t-1}{i-1}\binom{rj+t-1}{n-i-1}.\] 
\end{proof}

\begin{corollary}$(\cite{SA}, \ Lemma \ 5.1.)$ 
If $j=n-i-1,$ then $K-$ algebra 
$K[{\bf{\mathcal{A}}}]$ is Gorenstein ring. \\
\end{corollary}

Let $S$ be a standard graded $K-$algebra over a field $K$. Recall that
 the $a-$\emph{invariant} of $S,$
denoted $a(S),$ is the degree as a rational function of the Hilbert
 series of $S,$ see for instance
(\cite[p.  99]{V}). If $S$ is Cohen-Macaulay and $\omega_{S}$ is the
 canonical module of $S,$ then 
\[a(S)= {} - \min  \{i\ | \ (\omega_{S})_{i} \neq 0\},
\] 
see \cite[p.  141]{BH} and \cite[Proposition 4.2.3]{V}.
In our situation $S = K[{\bf{\mathcal{A}}}]$ is normal \cite{HH} and
consequently Cohen-Macaulay, thus this formula applies. We have the
 following consequence of Theorem 4.4.
\begin{corollary}
The $a-$invariant of $K[{\bf{\mathcal{A}}}]$ is
 $a(K[{\bf{\mathcal{A}}}])= \left \{ \begin{array}[pos]{cccc}
-1 \ , \  if \ \ \  i+j\leq n-1 \\

-r \ , \  if \ \ \  i+j\geq n \ \ \ \ \ 
\end{array} \right. $
\end{corollary}
\begin{proof}
Let $\{x^{\alpha_{1}}, \ldots, x^{\alpha_{q}}\}$ be  generators  of the
 $K-$algebra $K[{\bf{\mathcal{A}}}].$ Then  $K[{\bf{\mathcal{A}}}]$ is
 a standard graded algebra with the grading
\[K[{\bf{\mathcal{A}}}]_{i}=\sum_{|c|=i}K(x^{\alpha_{1}})^{c_{1}}\cdot
 \cdot \cdot (x^{\alpha_{q}})^{c_{q}}, \ 
\mbox{where} \ |c|=c_{1}+\ldots+c_{q}. 
\] 
If $i+j\leq n-1,$ then  
$ \min  \{i\ | \ (\omega_{K[{\bf{\mathcal{A}}}]})_{i} \neq 0\}=1$ and so
we have $a(K[{\bf{\mathcal{A}}}])=  -1.$\\
If $i+j\geq n,$ then  
$ \min  \{i\ | \ (\omega_{K[{\bf{\mathcal{A}}}]})_{i} \neq 0\}=r.$ 
So, we have
$a(K[{\bf{\mathcal{A}}}])=  -r.$
\end{proof}

\section{Ehrhart function}
We consider  a fixed set of distinct monomials
$F=\{x^{\alpha_{1}},\ldots,x^{\alpha_{r}}\}$ \ in a polynomial ring
$R=K[x_{1},\ldots, x_{n}]$ \ over a field $K.$ \\ 

Let
\[\mathcal{P}=conv(log(F))\] \ be the convex hull of the set 
$log(F)=\{\alpha_{1},\ldots,
\alpha_{r}\}.$ \ The $normalized \ Ehrhart \ ring$ \ of $\mathcal{P}$ is the
graded algebra \[A_{\mathcal{P}}=\bigoplus^{\infty}_{i=0}(A_{\mathcal{P}})_{i} 
\subset R[T]\] where the $i-th$ \ component  is given by
\[(A_{\mathcal{P}})_{i}=\sum_{\alpha \in \mathbb{Z} \ log(F)\cap \ 
i\mathcal{P}} \ K \ x^{\alpha} \ T^{i}.\]
The $normalized \ Ehrhart \ function$ of $\mathcal{P}$ is defined as
\[E_{\mathcal{P}}(t)= dim_{K}(A_{\mathcal{P}})_{t}=|\ \mathbb{Z} \ 
{log(F)}\cap \  t \ \mathcal{P} \ |.\] An important result of
{\cite{V}, Corollary 7.2.45} is the following :
\begin{theorem}
If $K[F]$ is a standard graded subalgebra of $R$ and $h$ \ is the
$Hilbert \ function$ of $K[F],$ \ then:\\
$a)$ \ $h(t)\leq E_{\mathcal{P}}(t)$ \ for all $t\geq 0,$ \ and\\
$b)$ \ $h(t)= E_{\mathcal{P}}(t)$ \ for all $t\geq 0$ \ if and only if $K[F]$
is normal.
\end{theorem}

In this section we will compute the Hilbert function and the Hilbert 
series for $K-$ algebra $K[{\bf{\mathcal{A}}}],$ where $\bf{\mathcal{A}}$ 
satisfied the hypothesis of $Lemma \ 4.1.$
\begin{proposition}
In the hypothesis of $Lemma \ 4.1.$, the Hilbert function of
 $K-$ algebra $K[{\bf{\mathcal{A}}}]$ is :
\[h(t)=\sum_{k=0}^{(n-j)t}\binom{k+i-1}{k}\binom{nt-k+n-i-1}{nt-k}.\]
\end{proposition}

\begin{proof}
From {\cite{HH}} we know that the $K-$ algebra $K[{\bf{\mathcal{A}}}]$ 
is normal. Thus, to compute the Hilbert function of $K[{\bf{\mathcal{A}}}]$ 
it is equivalent to compute the Ehrhart function of ${\mathcal{P}},$ 
where ${\mathcal{P}} \ = \ conv(A).$ $(see \ Theorem \ 5.1.)$ \\ 
It is clear that $\mathcal{P}=\{ x\in \mathbb{R}^n \ | \ 
x_{i}\geq 0, \ 0\leq x_{1}+\ldots + x_{i}\leq n-j \ and \ x_{1}+
\ldots + x_{n} = n \}$ and thus, $t \ \mathcal{P}=\{ x\in \mathbb{R}^n 
\ | \ x_{i}\geq 0, \ 0\leq x_{1}+\ldots + x_{i}\leq (n-j) \ t \ and 
\ x_{1}+\ldots + x_{n} = n \ t \}.$

Since for any $0\leq k \leq (n-j) \ t$ the equation $ x_{1}+\ldots + 
x_{i} = k$ has $\binom{k+i-1}{k}$ positive solutions and the equation 
$x_{i+1}+\ldots + x_{n} = n \ t -k$ has $\binom{nt-k+n-i-1}{nt-k}$ 
positive solutions, we get that \[E_{\mathcal{P}}(t)=|\ \mathbb{Z} 
\ A \ \cap \ t \ \mathcal{P} \ | = \sum_{k=0}^{(n-j)t}\binom{k+i-1}
{k}\binom{nt-k+n-i-1}{nt-k}.\]
\end{proof}

\begin{corollary}
In the hypothesis of $Lemma \ 4.1.$, the Hilbert series of 
$K-$ algebra $K[{\bf{\mathcal{A}}}]$ is :
\[H_{K[{\bf{\mathcal{A}}}]}(t) \ = \ \frac{1+h_{1} \ t + \ldots + 
h_{n-r} \ t^{n-r} }{(1-t)^{n}},\] where 
\[h_{j}=\sum_{s=0}^{j}(-1)^{s} \ h(j-s) \ \binom{n}{s},\]
$h(s)$ being the Hilbert function and $r=\left\lceil \frac{i+1}{n-j}\right\rceil.$
\end{corollary}

\begin{proof}
Since the $a-$invariant of $K[{\bf{\mathcal{A}}}]$ is
 $a(K[{\bf{\mathcal{A}}}])=  -r,$ with $r=\left\lceil 
 \frac{i+1}{n-j}\right\rceil,$ it follows that to 
 compute the Hilbert series of
 $K[{\bf{\mathcal{A}}}]$ it is necessary to know the first $n-r+1$ values of the
 Hilbert function of $K[{\bf{\mathcal{A}}}]$, $h(i)$ for $0\leq i \leq
 n-r.$
Since $\dim (K[{\bf{\mathcal{A}}}])=n,$ applying $n$ times the
 difference operator  $\Delta$ ({see  \cite{BH}})  on the Hilbert function of
 $K[{\bf{\mathcal{A}}}]$ we get the conclusion.

Let $\Delta^{0}(h)_{j}:=h(j)$ for any $0\leq j \leq n-r.$ For $k\geq 1$
 let $\Delta^{k}(h)_{0}:=1$
and $\Delta^{k}(h)_{j}:=\Delta^{k-1}(h)_{j}-\Delta^{k-1}(h)_{j-1}$
for any $1\leq j \leq n-r.$
We claim that
\[\Delta^{k}(h)_{j}=\sum_{s=0}^{k}(-1)^sh(j-s)\binom{k}{s}
\]
for any $k\geq 1$ and $0\leq j \leq n-r.$ We proceed by induction on
 $k.$

If $k=1,$ then
\[
\Delta^{1}(h)_{j}=\Delta^{0}(h)_{j}-\Delta^{0}(h)_{j-1}=h(j)-h(j-1)
=\sum_{s=0}^{1}(-1)^{s}h(j-s)\binom{1}{s}
\]
for any $1\leq j \leq n-r.$\\
If $k>1,$ then
\[\Delta^{k}(h)_{j}=\Delta^{k-1}(h)_{j}-\Delta^{k-1}
(h)_{j-1}=\sum_{s=0}^{k-1}(-1)^{s}h(j-s)
\binom{k-1}{s}-\sum_{s=0}^{k-1}(-1)^{s} h(j-1-s)
\binom{k-1}{s}\]
\[=h(j)\binom{k-1}{0}+\sum_{s=1}^{k-1}(-1)^{s}h(j-s)
\binom{k-1}{s}-\sum_{s=0}^{k-2}(-1)^{s} h(j-1-s)
\binom{k-1}{s}+(-1)^kh(j-k)\binom{k-1}{k-1}\]
\[=h(j)+\sum_{s=1}^{k-1}(-1)^{s}h(j-s)\left[\binom{k-1}{s}
+\binom{k-1}{s-1}\right]+(-1)^kh(j-k)\binom{k-1}{k-1}\
 \ \ \ \ \ \ \ \ \ \ \ \ \ \ \ \ \ \ \ \ \ \ \ \ \ \ \ \]
\[=h(j)+\sum_{s=1}^{k-1}(-1)^{s}h(j-s)\binom{k}{s}+(-1)^
kh(j-k)\binom{k-1}{k-1}=\sum_{s=0}^{k}(-1)^{s}h(j-s)\binom{k}{s}.\
 \ \ \ \ \ \ \ \ \ \ \ \ \ \ \ \ \ \ \ \ \ \ \ \]
Thus, if $k=n$ it follows that
\[h_{j}=\Delta^{n}(h)_{j}=\sum_{s=0}^{n}(-1)^sh(j-s)\binom{n}{s}=
\sum_{s=0}^{j}(-1)^sh(j-s)\binom{n}{s}\]
for any $1\leq j\leq n-r.$

\end{proof}
Next we will give some examples.
\begin{example}
Let ${\bf{\mathcal{A}}}=\{A_{1},\ldots ,A_{7}\},$
where $A_{1}=A_{2}=A_{3}=A_{6}=A_{7}=[7],$ \ $A_{4}=A_{5}=[7]\setminus [3].$
The cone genereated by $A,$ the exponent set of generators of
$K-$algebra $K[{\bf{\mathcal{A}}}],$ has the
irreducible representation:
\[\mathbb{R}_{+}A= H^{+}_{\nu^{2}_{\sigma^{0}[3]}}\cap 
H^{+}_{e_{1}} \cap \ldots \cap H^{+}_{e_{7}}.\] 

The type of $K[{\bf{\mathcal{A}}}]$ is:
\[type(K[{\bf{\mathcal{A}}}])=1+ \binom{8}
{2}\binom{4}{3}=113.\]

The Hilbert series of $K[{\bf{\mathcal{A}}}]$ is:
\[H_{K[{\bf{\mathcal{A}}}]}(t) \ = \ \frac{1+1561 t + 24795 t^{2} + 57023 t^{3}+
25571 t^{4} + 1673 t^{5} + t^{6}}{(1-t)^{7}}.\]
Note that
$type(K[{\bf{\mathcal{A}}}])=1+h_{5}-h_{1}=113.$
\end{example}

\begin{example}
Let ${\bf{\mathcal{A}}}=\{A_{1},\ldots ,A_{7}\},$
where $A_{3}=A_{4}=[7],$ \ $A_{1}=A_{2}=A_{5}=A_{6}=A_{7}=[7]\setminus [4].$
The cone genereated by $A,$ the exponent set of generators of
$K-$algebra $K[{\bf{\mathcal{A}}}],$ has the
irreducible representation:
\[\mathbb{R}_{+}A= H^{+}_{\nu^{5}_{\sigma^{0}[4]}}\cap 
H^{+}_{e_{1}} \cap \ldots \cap H^{+}_{e_{7}}.\] 

The type of $K[{\bf{\mathcal{A}}}]$ is:
\[type(K[{\bf{\mathcal{A}}}])=\binom{4}
{3}\binom{15}{2}+\binom{3}
{3}\binom{16}{2}=540.\]

The Hilbert series of $K[{\bf{\mathcal{A}}}]$ is:
\[H_{K[{\bf{\mathcal{A}}}]}(t) \ = \ \frac{1+351t + 2835t^{2} + 3297t^{3}+
540t^{4}}{(1-t)^{7}}.\]
\end{example}

Note that
$type(K[{\bf{\mathcal{A}}}])=h_{4}=540.$

We end this section with the following open problem:

{\bf{Open Problem}:}
Let $n\geq 4,$ $A_{i}\subset [n]$ for any $1\leq i \leq n$ and 
$K[{\bf{\mathcal{A}}}]$ be the base ring associated to 
the transversal polymatroid presented by
${\bf{\mathcal{A}}}=\{A_{1},\ldots ,A_{n}\}.$
If the Hilbert series is: 
\[H_{K[{\bf{\mathcal{A}}}]}(t) \ = \ \frac{1+h_{1} \ t + \ldots + 
h_{n-r} \ t^{n-r} }{(1-t)^{n}},\] then we have the followings:\\
$1)$ If $r=1,$ then $type(K[{\bf{\mathcal{A}}}])=1+h_{n-2}-h_{1}.$\\
$2)$ If $2\leq r\leq n,$ then $type(K[{\bf{\mathcal{A}}}])=h_{n-r}.$

{\bf Acknowledgments:}
The author whould like to thank Professor Dorin Popescu for valuable suggestions 
and comments during the preparation of this paper.

\vspace{2mm} \noindent {\footnotesize
\begin{minipage}[b]{10cm}
Alin \c{S}tefan, Assistant Professor\\
"Petroleum and Gas" University of Ploie\c{s}ti\\
Ploie\c{s}ti, Romania\\
E-mail:nastefan@upg-ploiesti.ro
\end{minipage}}
\end{document}